\documentclass{article}

\usepackage{mathptmx}
\usepackage{fancyhdr}
\usepackage{eucal}
\usepackage{amsmath}
\usepackage{amscd}
\usepackage{amssymb}
\usepackage{amsthm}
\usepackage{xspace}
\usepackage[all,tips]{xy}
\usepackage[dvips]{graphicx}
\usepackage{verbatim}
\usepackage{syntonly}

\theoremstyle{plain}
\newtheorem{thm}{Theorem}[section]

\newtheorem{cor}[thm]{Corollary}

\newtheorem{prop}[thm]{Proposition}
\newtheorem{lemma}[thm]{Lemma}

\theoremstyle{definition}
\newtheorem*{rem}{Remark}

\newenvironment{pf}
{\begin{proof}} {\end{proof}}

 \DeclareMathOperator{\Stab}{Stab}

\DeclareMathOperator*{\ch}{char}

\DeclareMathOperator{\GL}{GL} 
\DeclareMathOperator{\Ort}{O} \DeclareMathOperator{\PO}{PO}
 \DeclareMathOperator{\SO}{SO}

\DeclareMathOperator{\Nil}{Nil}

 \DeclareMathOperator{\Aff}{Aff}
\DeclareMathOperator{\Sim}{Sim} \DeclareMathOperator{\Euc}{Euc}

\DeclareMathOperator*{\Isom}{Isom}

\newcommand{\vp}{\varphi}

\newcommand{\al}{\alpha}

\newcommand{\be}{\beta}

\newcommand{\ga}{\gamma}
\newcommand{\Ga}{\Gamma}

\newcommand{\te}{\theta}

\newcommand{\la}{\lambda}
\newcommand{\La}{\Lambda}

\newcommand{\fn}{\footnote}

\newcommand{\iny}{\infty}

\newcommand{\tri}{\ensuremath{\triangle}}
\newcommand{\prt}{\partial}

\newcommand{\co}{\ensuremath{\colon}}

\newcommand{\innp}[1]{\left< #1 \right>}
\newcommand{\abs}[1]{\left\vert#1\right\vert}
\newcommand{\set}[1]{\left\{#1\right\}}

\newcommand{\su}{\subset}

\newcommand{\bu}{\bigcup}
\newcommand{\ba}{\bigcap}

\newcommand{\op}{\oplus}

\newcommand{\lra}{\longrightarrow}

\newcommand{\lmto}{\longmapsto}

\newcommand{\B}[1]{\ensuremath{\mathbf{#1}}}
\newcommand{\BB}[1]{\ensuremath{\mathbb{#1}}}
\newcommand{\Cal}[1]{\ensuremath{\mathcal{#1}}}
\newcommand{\Fr}[1]{\ensuremath{\mathfrak{#1}}}

\newcommand{\Hy}{\ensuremath{\B{H}}}

\newcommand{\Q}{\ensuremath{\B{Q}}}
\newcommand{\R}{\ensuremath{\B{R}}}
\newcommand{\Z}{\ensuremath{\B{Z}}}
\newcommand{\C}{\ensuremath{\B{C}}}
\newcommand{\RP}{\ensuremath{\B{RP}}}

\newcommand{\refP}[1]{Proposition~\ref{P:#1}}

\newcommand{\refT}[1]{Theorem~\ref{T:#1}}
\newcommand{\refL}[1]{Lemma~\ref{L:#1}}

\fancypagestyle{plain}

\begin{document}


\bibliographystyle{amsplain}


\title{\textbf{Arithmetic cusp shapes are dense}}
\author{D. B. McReynolds\fn{Supported by a V.I.G.R.E graduate fellowship and Continuing
Education fellowship.}} \maketitle


\begin{abstract}
\noindent In this article we verify an orbifold version of a
conjecture of Nimershiem from 1998. Namely, for every flat
$n$--manifold $M$, we show that the set of similarity classes of
flat metrics on $M$ which occur as a cusp cross-section of a
hyperbolic $(n+1)$--orbifold is dense in the space of similarity
classes of flat metrics on $M$. The set used for density is
precisely the set of those classes which arise in arithmetic
orbifolds.
 \end{abstract}

\section{Main results}

\noindent By a \emph{flat $n$--manifold} we mean a closed manifold
$M=\R^n/\Ga$ where $\Ga$ is a discrete, torsion free, cocompact
subgroup of $\Isom(\R^n)$. In analogy with Teichm\"{u}ller theory,
there is a contractible space $\Cal{F}(M)$ of flat metrics on $M$
coming from the standard flat structure on $\R^n$ and all the
possible $\Isom(\R^n)$--conjugacy classes for $\Ga$. We say that two
flat metrics $g_1,g_2$ on $M$ are \emph{similar} if there exists an
isometry between $(M,\al g_1)$ and $(M,g_2)$ for some $\al \in
\R^+$. We denote the equivalence class under similarity of a flat
metric $g$ by $[g]$ and the space of similarity classes of flat
metrics on $M$ by $\Cal{S}(M)$.\smallskip\smallskip

\noindent An important relationship between flat and hyperbolic
geometry is exhibited in the thick-thin decomposition of a
hyperbolic manifold. Specifically, every finite volume, noncompact
hyperbolic $(n+1)$--orbifold $W$ has a thick-thin decomposition
comprised of a compact manifold $W_{core}$ with boundary components
$M_1,\dots,M_m$ and manifolds $E_1,\dots,E_m$ of the form $M_j
\times \R^{\geq 0}$. The manifolds $E_j$ are called \emph{cusp
ends}, the manifolds $M_j$ are called \emph{cusp cross-sections},
and the union of $W_{core}$ along the boundary with $E_1,\dots,E_m$
recovers $W$ topologically. The manifolds $M_j$ are known to be flat
$n$--manifolds and are totally geodesic boundary components of the
manifold $W_{core}$ equipped with the quotient metric coming from
the path metric on its neutered space $\Cal{N} \su \Hy^{n+1}$.
Indeed, a cusp cross-section $M_j$ is naturally furnished with a
flat metric $g$ which is well-defined up to similarity and we call
these similarity classes of flat metrics \emph{realizable flat
similarity classes} or \emph{cusp shapes}. This article is devoted
to the classification of the possible cusp shapes of a flat
$n$--manifold occurring in the class of arithmetic
$(n+1)$--orbifolds. The absence of a general geometric construction
for hyperbolic orbifolds forces our restriction to orbifolds
produced by arithmetic means. Given this forced restriction, the
picture we provide here is complete. Before stating our main
results, we briefly survey some preexisting results and questions.
\smallskip\smallskip

\noindent Motivated by the above picture and work of Gromov
\cite{Gromov78}, Farrell and Zdravkovska \cite{FarrellZdravkovska83}
conjectured that every flat $n$--manifold arises as a cusp
cross-section of a 1--cusped hyperbolic $(n+1)$--manifold and this
is easily verified for $n=2$. Indeed, the complement of a knot in
$S^3$ is typically endowed with a finite volume, complete hyperbolic
structure with one cusp (see \cite{Thurston82}), and thus gives the
realization of the 2--torus $T^2$ as a cusp cross-section of a
1--cusped hyperbolic 3--manifold. Likewise, the Klein bottle arises
as a cusp cross-section of the 1--cusped Gieseking manifold (see
\cite{Ratcliffe94}). However, Long and Reid \cite{LongReid00}
constructed counterexamples for $n=3$ by showing that any flat
3--manifold arising as a cusp cross-section of a 1--cusped
hyperbolic 4--manifold must have integral $\eta$--invariant; this
works for all $n=4k-1$.\smallskip\smallskip

\noindent The failure of the conjecture of Farrell--Zdravkovska is
far from total. Nimershiem \cite{Nimershiem98} showed every flat
3--manifold arises as a cusp cross-section of a hyperbolic
4--manifold, and Long and Reid \cite{LongReid02} proved every flat
$n$--manifold arises as a cusp cross-section of an arithmetic
hyperbolic $(n+1)$--orbifold. A more geometrically relevant question
is precisely which cusp shapes occur on a given flat $n$--manifold.
Via a counting argument, almost every similarity class on a flat
$n$--manifold must fail to appear as a cusp shape. Despite this,
Nimershiem \cite{Nimershiem98} showed any for flat 3--manifold $M$
the cusp shapes occurring in hyperbolic 4--manifolds are dense in
the space of flat similarity classes on $M$, and conjectured
\cite[Conj. 2']{Nimershiem98} this for every flat $n$--manifold. Our
main result is the verification of this conjecture in the orbifold
category.

\begin{thm}[Cusp shape density]\label{T:SDensity}
For a flat $n$--manifold $M$, the set of cusp shapes of $M$
occurring in hyperbolic $(n+1)$--orbifolds is dense in the space of
flat similarity classes $\Cal{S}(M)$.
\end{thm}

\noindent It is worth mentioning that the proof of \refT{SDensity}
exhibits a dense subset of shapes of a uniform nature. These
similarity classes are the image of $\Q$--points of a
$\Q$--algebraic set under a projection map. From this one sees
density occurs not as a function of small complexity in low
dimensions but from the algebraic structure of these spaces.
Moreover, these similarity classes of flat metrics are precisely
those classes which occur in the cusp cross-sections of arithmetic
hyperbolic $(n+1)$--orbifolds---see
\refT{MetricClass}.\smallskip\smallskip

\noindent Using a modest refinement of Selberg's lemma, we verify
the full conjecture for the $n$--torus.

\begin{thm}\label{T:Density}
For the $n$--torus $\R^n/\Z^n$, the set of cusp shapes of
$\R^n/\Z^n$ occurring in hyperbolic $(n+1)$--manifolds is dense in
the space of flat similarity classes $\Cal{S}(\R^n/\Z^n)$.
\end{thm}

\subsection*{Acknowledgements}

\noindent I would like to thank my advisor Alan Reid for all his
help, and am grateful to Oliver Baues, Yves Benoist, Jim Davis,
Nathan Dunfield, Elisha Falbel, Yoshi Kamishima, Richard Kent, and
Chris Leininger for stimulating conversations on this work. Finally,
I thank Ernest Vinberg for pointing out some errors in an early
version of this article and the referee for numerous comments.

\section{Preliminaries}

\subsection{Bieberbach groups and flat manifolds}

\noindent We shall denote the affine, Euclidean, and similarity
groups of $\R^n$ by $\Aff(n)$, $\Euc(n)$, and $\Sim(n)$. We call a
discrete, torsion free subgroup $\Ga$ of $\Euc(n)$ a
\emph{Bieberbach group} when the quotient space $\R^n/\Ga$ is a
closed manifold and throughout the remainder of this article, $\Ga$
shall denote a Bieberbach group. We denote the space of faithful
representations $\rho$ of $\Ga$ in $\Aff(n)$ with Bieberbach images
by $\Cal{R}_f(\Ga)$ (i.e., $\rho(\Ga)$ is an $\Aff(n)$--conjugate of
a Bieberbach group $\Ga^\prime$). Finally, $\Cal{F}_f(\Ga),
\Cal{S}_f(\Ga)$ will denote the subspaces of the $\Euc(n)$ and
$\Sim(n)$--character spaces consisting of those faithful characters
whose images are Bieberbach.\smallskip\smallskip

\noindent Associated to each maximal compact subgroup $K$ of
$\GL(n;\R)$ is the \emph{orthogonal affine group} $\Ort_K(n)= \R^n
\rtimes K$. As each $K$ is conjugate in $\GL(n;\R)$ to $\Ort(n)$,
$K$ is equal to $\Ort(B_K)$ for some symmetric, positive definite,
bilinear form $B_K$ on $\R^n$. When $B_K$ is $\Q$--defined (i.e.
$\Q$--valued on some $\R$--basis), $\Ort_K(n)$ is $\Q$--defined and
we shall call subgroups commensurable with $\Ort_K(n;\Z)$
\emph{$\Q$--arithmetic subgroups} of $\Ort_K(n)$. We say $\rho$ in
$\Cal{R}_f(\Ga)$ is \emph{$\Q$--arithmetic} if there exists a
$\Q$--defined orthogonal affine group $\Ort_K(n)$ such that
$\rho(\Ga)$ is a $\Q$--arithmetic subgroup of $\Ort_K(n)$, and
denote the subspace of $\Cal{R}_f(\Ga)$ of $\Q$--arithmetic
representations by $\Cal{R}_f(\Ga;\Q)$.\smallskip\smallskip

\noindent Recall that a manifold $M$ is \emph{flat} if $M$ is
diffeomorphic to $\R^n/\Ga$ for some Bieberbach group $\Ga$. The
flat metric $g$ induced by the standard inner product
$\innp{\cdot,\cdot}$ on $\R^n$ supplies $M$ with a flat metric
called the \emph{associated flat structure}. For future reference
$\Cal{F}(M)$ shall denote the space of all isometry classes of such
metrics.

\begin{thm}[\cite{Thurston97}]
The space of flat isometry classes on $M$ is $\Cal{F}_f(\pi_1(M))$.
The space of flat similarity classes on $M$ is
$\Cal{S}_f(\pi_1(M))$.
\end{thm}

\noindent These identifications are as real analytic spaces. Also,
observe that any faithful representation $\rho$ of $\pi_1(M)$ into
$\Ort_K(n)$ whose image is Bieberbach endows $M$ with a flat metric
induced from the form $B_K$.

\subsection{Hyperbolic space $\Hy^{n+1}$}

\noindent The classical group $\SO(n+1,1)$ produces the symmetric
space $\Hy^{n+1}$ known as hyperbolic $(n+1)$--space. For an
explicit description, we equip $\R^{n+2}$ with a bilinear form $B$
of signature $(n+1,1)$. \emph{Hyperbolic $(n+1)$--space} is the
$\R$--projectivization of the space of $B$--negative vectors (i.e.,
$v \in \R^{n+2}$ such that $B(v,v)<0$) endowed with the Bergman
metric associated to $B$. We denote hyperbolic $(n+1)$--space
together with this metric by $\Hy^{n+1}$ and say $\Hy^{n+1}$ is
\emph{modelled on $B$} and call $B$ a \emph{model form}. The
\emph{boundary} of $\Hy^{n+1}$ in $\RP^{n+1}$ is the
$\R$--projectivization of the space of $B$--null vectors (i.e., $v
\in \R^{n+2}$ such that $B(v,v)=0$) and we denote this set by $\prt
\Hy^{n+1}$.\smallskip\smallskip

\noindent The isometry group of hyperbolic $(n+1)$--space
$\Hy^{n+1}$ is afforded a $KAN$ decomposition called \emph{the
Iwasawa decomposition}. The decomposition depends on the selection
of a pair of boundary points $[v_0],[v_\iny]$. Most important for us
here is the factor $N$ which is isomorphic to $\R^n$ (the factor $K$
is a maximal compact subgroup and the factor $A$ consists of all
hyperbolic transformations fixing $[v_0]$ and $[v_\iny]$). Briefly,
all such isomorphisms arise as follows. Let $B$ be a model bilinear
form for hyperbolic $(n+1)$--space and $V_\iny$ be the
$B$--orthogonal complement of a pair of linearly independent
$B$--null vectors $v_0$ and $v_\iny$ in $\R^{n+2}$. For any maximal
compact $K$ of $\GL(n;\R)$ with associated positive definite form
$B_K$, let
\[ \psi\co (\R^n,B_K) \lra (V_\iny,B_{|V_\iny}) \]
be any isometric isomorphism. This induces the desired isomorphism
$\eta\co \R^n \lra N$ defined by
\[ \eta(\xi) = \exp(\psi(\xi) v_\iny^* - v_\iny\psi(\xi)^*), \]
where $xy^*(\cdot) = B(\cdot,y)x$ is the outer pairing of $x$ and
$y$ with respect to the form $B$. The injectivity and linearity of
this map are straightforward to check. That this is surjective
(i.e., $N$ is isomorphic to $\R^n$) is less clear but well known and
can be found in \cite{Ratcliffe94}. This extends to $\eta\co
\Ort_K(n) \lra \Isom(\Hy^{n+1})$, and produces
\[ \eta(\Ort_K(n-1)) < \Stab([v_\iny]). \]\smallskip\smallskip

\noindent Viewing $\Hy^{n+1}$ as the coset space
$\Isom(\Hy^{n+1})/K$, where $K$ is a maximal compact subgroup,
$\Hy^{n+1}$ is identified with $A \times N$ (the latter can be given
a natural $\Isom(\Hy^{n+1})$--invariant metric for which this
identification is an isometry). As $A=\R^+$, $\Hy^{n+1}$ supports
the foliation
\[ \bu_{t \in \R^+} \set{t} \times N \]
whose leaves are called \emph{horospheres} and are said to be
\emph{centered at} $[v]$ if $[v] \in \prt \Hy^{n+1}$ plays the role
of $v_\iny$ in the discussion above on the isomorphism between $N$
and $\R^n$. \smallskip\smallskip

\noindent By a \emph{lattice} in $\Isom(\Hy^{n+1})$, we mean a
discrete subgroup $\La$ such that $\Hy^{n+1}/\La$ has finite volume.
If $\Hy^{n+1}/\La$ is compact, we say that $\La$ is \emph{cocompact}
and otherwise call $\La$ \emph{noncocompact}. We call the quotient
$\Hy^{n+1}/\La$ a \emph{hyperbolic $(n+1)$--orbifold} and when $\La$
is arithmetic, the quotient is referred to as an \emph{arithmetic
hyperbolic $(n+1)$--orbifold}. Of primary concern in this article
are arithmetic lattices whose description is given by the following
well known theorem---see for instance \cite{LiMillson93} or
\cite{McRThesis}. For the statement, recall that a pair of subgroups
$H_1,H_2$ of a group $G$ are \emph{commensurable in the wide sense}
if there exists $g \in G$ such that $(g^{-1}H_1g) \cap H_2$ is a
finite index subgroup of both $g^{-1}H_1g$ and $H_2$ (when $g$ is
trivial, we say $H_1$ and $H_2$ are commensurable).

\begin{thm}\label{T:Arithmetic}
Every noncocompact arithmetic lattice in $\Isom(\Hy^{n+1})$ is
commensurable in the wide sense with $\PO_0(B;\Z)$ for some
$\Q$--defined signature $(n+1,1)$ form $B$. Conversely, for $n \geq
3$, any such $\PO_0(B;\Z)$ is a noncompact arithmetic lattice in
$\Isom(\Hy^{n+1})$.
\end{thm}

\noindent For a lattice $\La$ in $\Isom(\Hy^{n+1})$ with associated
orbifold $W=\Hy^{n+1}/\La$, we say that $W$ has a \emph{cusp} at
$[v]$ if $\La \cap N \ne \set{1}$, where $N$ is the group associated
to $[v]$. In this case, the \emph{maximal peripheral subgroup of
$\La$ at $[v]$} is the subgroup $\tri_v(\La) = \Stab([v]) \cap \La$.
By the Kazhdan--Margulis theorem, $\tri_v(\La)$ is virtually abelian
and a maximal abelian subgroup of $\tri_v(\La)$ is given by $L =
\tri_v(\La) \cap N$. Moreover, the Kazhdan--Margulis theorem permits
the selection of a horosphere $\set{t}\times N$ such that
$\set{t}\times N/\tri_v(\La)$ is embedded in $\Hy^{n+1}/\La$. In
this case, we call $\set{t}\times N/\tri_v(\La)$ a \emph{cusp
cross-section} of the cusp at $[v]$.\smallskip\smallskip

\noindent Each cusp cross-section is equipped with a flat structure
coming from the orthogonal affine group $\Ort_K(n)$ of
$\Isom(\Hy^{n+1})$ containing $\tri_v(\La)$. If $B$ is a signature
$(n+1,1)$ model form for $\Hy^n$, the associated form producing
$\Ort_K(n)$ is obtained by selecting a different boundary point
$[w]$ and considering the restriction of $B$ to the $B$--orthogonal
complement of the $\R$--span of $v$ and $w$. Up to similarity, this
flat structure is independent of our selection of $[w]$ and thus
yields a point in $\Cal{S}(M)$. We say that $[g] \in \Cal{S}(M)$ is
\emph{realizable} or is a \emph{cusp shape} if there exists a
hyperbolic $(n+1)$--orbifold $W$ with cusp cross-section
$(M^\prime,[g^\prime])$ and a similarity transformation $f\co
(M,[g]) \lra (M^\prime,[g^\prime])$.

\section{Proofs of \refT{SDensity} and \refT{Density}}

\subsection{Density of shapes}

\noindent The proof of \refT{SDensity} is established in three
steps. We first show that each representation $\rho$ in
$\Cal{R}_f(\Ga;\Q)$ produces a similarity class that arises as a
cusp shape in an arithmetic hyperbolic $(n+1)$--orbifold. After
verifying the density of $\Cal{R}_f(\Ga;\Q)$ in $\Cal{R}_f(\Ga)$, we
construct continuous surjective maps from $\Cal{R}_f(\Ga)$ to both
$\Cal{F}_f(\Ga)$ and $\Cal{S}_f(\Ga)$. \smallskip\smallskip

\begin{prop}\label{P:StrongClass}
Let $\rho \in \Cal{R}_f(\Ga;\Q)$ and $\Ort_K(n)$ be any
$\Q$--defined orthogonal affine group with $\rho(\Ga)$ a
$\Q$--arithmetic subgroup of $\Ort_K(n)$. Then there exists an
arithmetic lattice $\La$ in $\Isom(\Hy^{n+1})$ and an injection
$\psi\co \Ort_K(n) \lra \Isom(\Hy^{n+1})$ such that
$\psi(\rho(\Ga))$ is a maximal peripheral subgroup of $\La$.
\end{prop}

\noindent Note that by Mal'cev rigidity, $\psi$ induces a similarity
transformation between $\R^n/\Ga$ and $\set{t}\times
N/\psi(\rho(\Ga))$. Consequently, in the sequel the transition from
algebraic statements to geometric ones will be made without comment.

\begin{pf}
We commence the proof by noting that with the peripheral
separability theorem in \cite{McReynolds04A}, we require only that
$\psi(\rho(\Ga))$ reside in $\La$. To see this, assume that $\psi$
is an injection of $\rho(\Ga)$ into an arithmetic lattice $\La_0$
and $\tri$ is the maximal peripheral subgroup containing
$\psi(\rho(\Ga))$. As $\psi(\rho(\Ga))$ is a finite index subgroup
of $\tri$, there exists a complete set of coset representatives
$\ga_1,\dots,\ga_r$ in $\tri$ for $\tri/\psi(\rho(\Ga))$. For each
$\ga_j$, the peripheral separability theorem of \cite{McReynolds04A}
provides us with a finite index subgroup $\La_j$ of $\La_0$ such
that $\psi(\rho(\Ga))$ is contained in $\La_j$ and $\ga_j \notin
\La_j$. Setting
\[ \La = \ba_j \La_j, \]
we obtain the needed pair $\psi$ and $\La$ for the validity of the
proposition. \smallskip\smallskip

\noindent It remains to find $\La_0$ and $\psi$, a task achieved
with the following line of reasoning. For the orthogonal affine
group $\Ort_K(n)$ with form $B_K$, define a model form $B$ for
$\Hy^{n+1}$ by
\[ B = B_K \op \begin{pmatrix} 1 & 0 \\ 0 & -1 \end{pmatrix} \]
and note that $B$ is $\Q$--defined being a direct sum of such forms.
We note that there is a $\Q$--defined injection of $\Ort_K(n)$ into
$\PO_0(B)$ defined by sending $\R^n$ to the group $N$ associated to
the boundary point $[e_{n+1}+e_{n+2}]$ and $K$ to the orthogonal
group acting trivially on the 2--plane spanned by $e_{n+1},e_{n+2}$.
The former is given explicitly by
\[ v \lmto \exp(v(e_{n+1}+e_{n+2})^* - (e_{n+1}+e_{n+2})v^*), \]
where we identify $\R^n$ inside $\R^{n+2}$ by taking the $n$--plane
spanned by $e_1,\dots,e_n$. The lattice $\La_0$ in $\PO_0(B)$ is
obtained from an application of \cite[Cor. 10.14]{Raghunathan72}.
Specifically, there exists $\La_0$ in $\PO_0(B)$ commensurable with
$\PO_0(B;\Z)$ such that $\psi(\rho(\Ga))$ is contained in $\La_0$.
\end{pf}

\begin{prop}\label{P:Orbit}
For every Bieberbach group $\Ga$, $\Cal{R}_f(\Ga;\Q)$ is dense in
$\Cal{R}_f(\Ga)$.
\end{prop}

\noindent To prove \refP{Orbit}, we need a pair of auxiliary lemmas.

\begin{lemma}\label{L:DenseOrbit}
If $G$ is a topological group with dense subgroup $H$ and $X$ is a
topological space with a continuous transitive $G$--action, then the
$H$--orbit of any $x$ in $X$ is dense in $X$.
\end{lemma}

\begin{lemma}\label{L:QPoints}
For every crystallographic group $\Ga$, $\Cal{R}_f(\Ga;\Q)$ is
nonempty.
\end{lemma}

\noindent The first lemma is elementary and easily verified. The
second lemma is a consequence of the Bieberbach theorems. Indeed, by
the Bieberbach theorems, one is provided with an injection
\[ \vp\co \Ga \lra \Z^n \rtimes \GL(n;\Z). \]
One can then take the $\te$--average (see below for more on this) of
the standard positive definite form on $\R^n$ to obtain a
$\Q$--defined positive definite form $B$ for which $\vp(\Ga)$ is an
arithmetic subgroup of $\Ort_{\Ort(B)}(n)$.

\begin{pf}[Proof of \refP{Orbit}]
As expected, we seek to apply \refL{DenseOrbit} and must ensure that
the conditions are satisfied by $X=\Cal{R}_f(\Ga)$ and $G=\Aff(n)$.
To begin, the topology on $\Cal{R}_f(\Ga)$ is the subspace topology
induced by viewing $\Cal{R}_f(\Ga)$ as a subspace of the
$\Aff(n)$--representation space. Visibly, the $\Aff(n)$--action on
the representation space is continuous, and so by restriction the
action of $\Aff(n)$ on $\Cal{R}_f(\Ga)$ is continuous. Less obvious
is the transitivity of the $\Aff(n)$--action on $\Cal{R}_f(\Ga)$.
However, this is precisely the statement of one part of the
Bieberbach theorems. Thus, by \refL{DenseOrbit}, for $H=\Q^n \rtimes
\GL(n;\Q)$ and any $\rho$ in $\Cal{R}_f(\Ga)$, the $H$--orbit of
$\rho$ is dense in $\Cal{R}_f(\Ga)$. We assert that for each $\al$
in $H$ and $\rho \in \Cal{R}_f(\Ga;\Q)$, the $\al$--conjugate
representation $\mu_\al \circ \rho$ is in $\Cal{R}_f(\Ga;\Q)$, where
\[ \mu_\al\co \Aff(n) \lra \Aff(n) \]
is the inner automorphism given by conjugating by $\al$. To see
this, let $\Ort_K(n)$ be a $\Q$--defined orthogonal affine group for
which $\rho(\Ga)$ is a $\Q$--arithmetic subgroup of $\Ort_K(n)$.
Conjugation by $\al$ yields an isomorphism
\[ \mu_\al\co \textrm{O}_K(n) \lra \textrm{O}_{\be^{-1}K\be}(n) \]
where $\be \in \GL(n;\Q)$ is the linear factor (or second
coordinate) for $\al$. As $\be$ is an element of $\GL(n;\Q)$, the
symmetric, positive definite form associated to $\be^{-1}K\be$ is
$\Q$--defined, being $\Q$--equivalent to the $\Q$--defined form
$B_K$. Moreover, this isomorphism between $\Ort_K(n)$ and
$\Ort_{\be^{-1}K\be}(n)$ is $\Q$--defined. Therefore, any
$\Q$--arithmetic subgroup of $\Ort_K(n)$ is mapped to a
$\Q$--arithmetic subgroup of $\Ort_{\be^{-1}K\be}(n)$, and thus
$\mu_\al \circ \rho(\Ga)$ is a $\Q$--arithmetic subgroup of a
$\Q$--defined orthogonal affine group as asserted. By
\refL{QPoints}, $\Cal{R}_f(\Ga;\Q)$ is nonempty, and so there exists
a $\Cal{R}_f(\Ga)$--dense $H$--orbit of representations in
$\Cal{R}_f(\Ga;\Q)$.
\end{pf}

\begin{pf}[Proof of \refT{SDensity}]
To begin, there exists a continuous surjective map
\[ \Cal{L}\co \Cal{R}_f(\Ga) \lra \Cal{F}_f(\Ga) \]
given as follows. For $\rho$ in $\Cal{R}_f(\Ga)$, as $\rho(\Ga)$ is
an $\Aff(n)$--conjugate of a Bieberbach group, by the Bieberbach
theorems $\rho(\Ga)$ projects to a finite group $\te$ in
$\GL(n;\R)$. Taking the $\te$--average
\[ B_\te(x,y) = \frac{1}{\abs{\te}}\sum_{g \in \te} \innp{gx,gy} \]
of the standard inner product $\innp{\cdot,\cdot}$ on $\R^n$
produces a maximal compact subgroup $K=\Ort(B_\te)$ such that
$\rho(\Ga)$ is contained in the orthogonal affine group $\Ort_K(n)$.
Up to post-composition with inner automorphisms of $\Euc(n)$, there
exists a unique $S_\rho$ in $\GL(n;\R)$ conjugating $\Ort_K(n)$ to
$\Euc(n)$. From this, we define $\Cal{L}(\rho) = S^{-1}_\rho \rho
S_\rho$.
\smallskip\smallskip

\noindent If $\rho_n$ is a sequence of representations in
$\Cal{R}_f(\Ga)$ converging to $\rho$ in $\Cal{R}_f(\Ga)$, the
sequence $\Cal{L}(\rho_n)$ converges to  $\Cal{L}(\rho)$ in
$\Cal{F}_f(\Ga)$. For a free abelian group this is immediate since
$\te$ is trivial. For nontrivial $\te$, this follows from the
convergence of the maximal compact subgroups $K_n$ arising from the
$\te_n$--average. We briefly explain this. For a sequence of
representations $\rho_n$ converging to $\rho$, let $\te_n$ be the
image of $\rho_n(\Ga)$ under projection onto $\GL(n;\R)$. It follows
that the groups $\te_n$ converge to $\te$ and thus $B_{\te_n}$
converges to $B_\te$ in the space of positive definite, symmetric
matrices---we use the standard basis to associate these matrices to
the $\te$--average forms. From this we see that the maximal compact
subgroups $K_n$ converge to $K$. The conjugating matrices
$S_{\rho_n}$ need not converge to the conjugating matrix $S_\rho$.
However, up to left multiplication in $\Ort(n)$, the sequence does
converge and so the sequence of $S_{\rho_n}^{-1}\rho_nS_{\rho_n}$
converges to $S_\rho^{-1}\rho S_\rho$ in $\Cal{F}_f(\Ga)$. As this
is the sequence $\Cal{L}(\rho_n)$, we see that $\Cal{L}(\rho_n)$
converges to $\Cal{L}(\rho)$ in $\Cal{F}_n(\Ga)$.
\smallskip\smallskip

\noindent The desired set of flat similarity classes is the image of
$\Cal{L}(\Cal{R}_f(\Ga;\Q))$ under the projection map
\[ \textrm{Pr}\co \Cal{F}_f(\Ga) \lra \Cal{S}_f(\Ga). \]
That this subset is dense is a consequence of the continuity and
surjectivity of $\Cal{L}$ in combination with \refP{Orbit}. For the
latter, if $\rho(\Ga)$ resides in $\Euc(n)$, the $\te$--average of
the standard form is the standard form and thus produces $\Ort(n)$.
In particular, we can take $S_\rho = I_n$ and hence $\Cal{L}$
restricted to subset of $\Euc(n)$--representations is the standard
projection onto $\Cal{F}_f(\Ga)$. It remains to show that each
similarity class in $\textrm{Pr}(\Cal{L}(\Cal{R}_f(\Ga;\Q)))$ does
occur as a cusp shape of an arithmetic hyperbolic $(n+1)$--orbifold.
By \refP{StrongClass}, for each $\rho \in \Cal{R}(\Ga;\Q)$ with
associated $\Q$--defined orthogonal affine group $\Ort_K(n)$, there
exists a faithful representation
\[ \psi\co \Ort_K(n) \lra \Isom(\Hy^{n+1}) \]
and an arithmetic lattice $\La$ such that $\psi(\rho(\Ga))$ is a
maximal peripheral subgroup of $\La$. For the flat structure on
$\R^n/\Ga$ coming from $\Ort_K(n)$ and the flat structure on the
cusp cross-section associated to $\psi(\rho(\Ga))$, this produces a
similarity of this pair of flat manifolds. To obtain this for the
associated class in $\textrm{Pr}(\Cal{L}(\rho))$, we argue as
follows. It could be that the $\Q$--form $B_K$ for $K$ is not the
$\te$--average of $\rho(\Ga)$. If this is the case, simply replace
$K$ by $\Ort(B_\te)$, and notice that this too is a $\Q$--defined
orthogonal affine group for which $\rho(\Ga)$ is a $\Q$--arithmetic
subgroup. Let $M^\prime$ be the associated flat manifold with this
similarity class associated to $\rho$ viewed as a representation
into $\Ort_{\Ort(B_\te)}(n)$. Making the same argument as before, we
see that $M^\prime$ occurs as a cusp shape of an arithmetic
hyperbolic $(n+1)$--orbifold. By construction, the flat manifold
$M^{\prime\prime}$ with similarity class
$\textrm{Pr}(\Cal{L}(\rho))$ is similar to $M^\prime$. Hence, every
class in the dense subset $\textrm{Pr}(\Cal{L}(\Cal{R}_f(\Ga;\Q)))$
arises as a cusp shape of an arithmetic real hyperbolic
$(n+1)$--orbifold.
\end{pf}

\begin{rem}
Yves Benoist pointed out to us that by using only the subgroup
separability result of \cite{McReynolds04A}, one can achieve density
in a fixed commensurability class of arithmetic hyperbolic
$(n+1)$--orbifolds for the $n$--torus. This argument can be extended
to other flat manifolds given the work of Long--Reid
\cite{LongReid02}.
\end{rem}

\noindent For $X = \C$ or $\BB{H}$, cusp cross-sections of finite
volume $X$--hyperbolic $(n+1)$--orbifolds are almost flat orbifolds
modelled on the $(2n+1)$--dimensional Heisenberg group
$\Fr{N}_{2n+1}$ or the $(4n+3)$--dimensional quaternionic Heisenberg
group $\Fr{N}_{4n+3}(\BB{H})$. Generalizing the approach of Long and
Reid, \cite{McReynolds04A} gave the smooth classification of cusp
cross-sections of $X$--hyperbolic $n$--orbifolds. Provided that an
almost flat manifold arises topologically as a cusp cross-section,
density follows with essentially the same argument.

\begin{thm}\label{T:AlmostFlat}
\begin{description}
\item[(a)]
For an almost flat $(2n-1)$--manifold $N$ modelled on
$\Fr{N}_{2n-1}$, the space of realizable almost flat similarity
classes in the cusp cross-sections of arithmetic complex hyperbolic
$n$--orbifolds is either empty or dense in the space of almost flat
similarity classes.
\item[(b)]
For an almost flat $(4n-1)$--manifold $N$ modelled on
$\Fr{N}_{4n-1}(\BB{H})$, the space of realizable almost flat
similarity classes in the cusp cross-sections of quaternionic
hyperbolic $n$--orbifolds is either empty or dense in the space of
almost flat similarity classes.
\end{description}
\end{thm}

\noindent As every $\Nil$ 3--manifold is diffeomorphic to a cusp
cross-section of an arithmetic complex hyperbolic $2$--orbifold (see
\cite{McReynolds04A}), \refT{AlmostFlat} yields:

\begin{cor}
For a $\Nil$ 3--manifold $N$, the space of $\Nil$ similarity classes
that arise in the cusp cross-sections of arithmetic complex
hyperbolic $2$--orbifolds is dense in the space of $\Nil$ similarity
classes.
\end{cor}

\noindent Finally, with the easily established converse of
\refP{StrongClass}, we obtain our main theorem, the geometric
classification of cusp cross-sections of arithmetic hyperbolic
$(n+1)$--orbifolds.

\begin{thm}[Geometric classification theorem]\label{T:MetricClass}
For a flat $n$--manifold $M$, the set \newline
$\textrm{Pr}(\Cal{L}(\Cal{R}_f(\pi_1(M);\Q)))$ is precisely the set
of flat similarity classes on $M$ that arise in cusp cross-sections
of arithmetic hyperbolic $(n+1)$--orbifolds.
\end{thm}

\noindent This persists in the complex and quaternionic hyperbolic
settings upon taking into account the above dichotomy.

\subsection{Orbifold to manifold promotion}

\noindent \refT{Density} is a consequence of the following
proposition whose proof is essentially a reproduction of Borel's
proof of Selberg's lemma \cite[Prop. 2.2]{Borel63}.

\begin{prop}\label{P:UnipotentSelberg}
If $k/\Q$ is a finite extension and $\La$ a finitely generated
subgroup of $\GL(n;k)$ with unipotent subgroup $\Ga$, then there
exists a torsion free, finite index subgroup $\La_0$ of $\La$ such
that $\Ga$ is contained in $\La_0$.
\end{prop}

\begin{pf}
Let $\la_1,\dots,\la_r$ be a finite generating set for $\La$,
$c_{i,j,\ell}$ be the $(i,j)$--coefficient of $\la_\ell$, and $R$ be
the subring of $k$ generated by $\set{c_{i,j,\ell}}$. By assumption,
$\Ga$ is conjugate in $\GL(n;\C)$ into the group of upper triangular
matrices with ones along the diagonal. In particular, the
characteristic polynomial $p_\ga(t)$ for each $\ga$ in $\Ga$ is
$(t-1)^n$. For any torsion element $\eta$ in $\La$, the
characteristic polynomial $p_\eta(t)$ of $\eta$ has only roots of
unity for its zeroes. Since $k/\Q$ is a finite extension and $n$ is
fixed, there are only finitely many degree $n$ monic polynomials in
$k[t]$ having only roots of unity for their roots. Let
$p_1(t),\dots,p_s(t)$ denote those monic polynomials with
coefficients in $R$ with this property. As our concern is solely
with nontrivial torsion elements, we further insist that each of the
polynomials has a root distinct from $1$. For each such polynomial
$p_j(t)$, there are finitely many prime ideals $\Fr{p}$ of $R$ such
that $(t-1)^n = p_j(t)$ modulo $\Fr{p}$. To see this, we first
exclude all prime ideals $\Fr{p}$ in $R$ such that
$\ch(R/\Fr{p})\leq n$. Since for each prime $p$ of $\Z$ there are
finitely many prime ideals $\Fr{p}$ of $R$ such that
$\ch(R/\Fr{p})=p$, this is a finite set. Next, as
\[ p_j(t) - (t-1)^n = \sum_{m=1}^n \la_{m,j}t^m \]
is nonzero, there exists $i$ such that $\la_{i,j}$ is nonzero. Since
$\Cal{O}_k$ is Dedekind, there are only finitely many prime ideals
$\Fr{p}_{j,1},\dots,\Fr{p}_{j,\ell_j}$ such that $\la_{i,j} = 0 \mod
\Fr{p}_{j,\ell_m}$ ($m=1,\dots,\ell_j$), and so for any other prime
ideal $\Fr{q}$, it follows that $p_j(t) \ne (t-1)^n$ modulo
$\Fr{q}$. Excluding this finite collection $\Cal{P}_j$ of prime
ideals of $R$, for any selection $\Fr{q} \notin \Cal{P}_j$, we have
$p_j(t)$ not equal to $(t-1)^n$ modulo $\Fr{q}$. Repeating this
argument for each $j$, we obtain the desired ideal set $\Cal{P}$.
For $\Fr{q} \notin \Cal{P}$, consider the reduction map $r_\Fr{q}\co
\GL(n;R) \lra \GL(n;R/\Fr{q})$. By our selection of $\Fr{q}$, no
torsion element $r_\Fr{q}(\eta)$ can be contained in $r_\Fr{q}(\Ga)$
as every element of $r_\Fr{q}(\Ga)$ has characteristic polynomial
$(t-1)^n$ and $\eta$ does not share this trait. Therefore,
$r_\Fr{q}^{-1}(r_\Fr{q}(\Ga))$ is a torsion free finite index
subgroup of $\La$ containing $\Ga$, as sought.
\end{pf}

\begin{rem}
We note that by Weil's local rigidity theorem and the Zariski
density of algebraic points, the assumption $\La<\GL(n;k)$ for a
finite extension $k/\Q$ is unnecessary.
\end{rem}

\begin{pf}[Proof of \refT{Density}]
It suffices to show that for each $\rho$ in $\Cal{R}(\Z^n;\Q)$ the
induced representation given by \refP{StrongClass} is such that
$\Z^n$ is contained in a torsion free finite index subgroup of the
target lattice $\La$. By construction, the representation
$\rho\co\Z^n \lra \La$ maps $\Z^n$ into a unipotent subgroup of
$\La$ since the groups $N$ in the Iwasawa decomposition are
unipotent. As the target lattice $\La$ is arithmetic, $\La$ is
finitely presentable (\cite[Cor. 13.25]{Raghunathan72}) and
conjugate into the $k$--points of $\Isom(\Hy^{n+1})$ for some number
field $k$. Thus \refP{UnipotentSelberg} is applicable and yields a
torsion free finite index subgroup $\La_0$ of $\La$ such that $\Z^n$
is contained in $\La_0$. Note that if $\Z^n$ is a maximal peripheral
subgroup of $\La$, $\Z^n$ is a maximal peripheral subgroup of
$\La_0$. In particular, we can realize the associated flat
similarity class $[g]$ for $\rho$ in a cusp cross-section of the
associated arithmetic manifold for $\La_0$.
\end{pf}

\noindent For those infranil manifold groups realizable as lattices
in their associated nilpotent Lie group, we say that the associated
infranil manifold is a \emph{niltorus}. For niltori modelled on
either $\Fr{N}_{2n-1}$ or $\Fr{N}_{4n-1}(\BB{H})$, orbifold density
is promoted to manifold density with an identical
argument.\smallskip\smallskip


\def\cprime{$'$} \def\lfhook#1{\setbox0=\hbox{#1}{\ooalign{\hidewidth
  \lower1.5ex\hbox{'}\hidewidth\crcr\unhbox0}}} \def\cprime{$'$}
  \def\cprime{$'$}
\providecommand{\bysame}{\leavevmode\hbox
to3em{\hrulefill}\thinspace}
\providecommand{\MR}{\relax\ifhmode\unskip\space\fi MR }
\providecommand{\MRhref}[2]{%
  \href{http://www.ams.org/mathscinet-getitem?mr=#1}{#2}
} \providecommand{\href}[2]{#2}


\noindent
Department of Mathematics, \\
California Institute of Technology\\
Pasadena, CA 91125\\
email: {\tt dmcreyn@caltech.edu}



\end{document}